\documentclass[10pt]{article}
\usepackage{amssymb,amsmath,amsthm,latexsym}
\usepackage{amsfonts}
\usepackage[twlrm]{rawfonts}
\usepackage[english]{babel}
\usepackage[latin1]{inputenc}
\usepackage{indentfirst}
\usepackage{color}
\newtheorem{teo}{Theorem}[section]
\newtheorem{cor}{Corolary}[section]
\newtheorem{lem}{Lemma}[section]

\newenvironment{dem}[1][Proof]{\noindent\textbf{#1.} }{\hfill \rule{0.5em}{0.5em}}

\newcommand{\N}{\mathbb{N}}

\newcommand{\R}{\mathbb{R}}

\newtheorem{ob}{Remark}[section]

\newtheorem{Af}{Claim}[section]

\begin{document}

%\chead{\thepage}

%\fancyhead[CO]{{\footnotesize \author{\rightmark}}\hspace{2em}\thepage}
%\setcounter{tocdepth}{2}
%\fancyhead[CE]{\thepage\hspace{2em}\footnotesize{\leftmark}}
%\fancyhead[CE,CO]{}
%\fancyhead[RO]{{\footnotesize\rightmark}\hspace{2em}\thepage}

%\setcounter{tocdepth}{2}\%fancyhead[LE]{\thepage\hspace{2em}\footnotesize{\leftmark}}
%\fancyhead[RE,LO]{}
%\fancyhead[RO]{{\footnotesize\rightmark}\hspace{2em}\thepage}	
	
\setlength{\baselineskip}{6.5mm} \setlength{\oddsidemargin}{8mm}
\setlength{\topmargin}{-3mm}

\title{\bf Existence of positive solution for a system of elliptic equations via bifurcation theory}

\author{Romildo N. de Lima\thanks{R. N. de Lima was partially supported by CAPES/Brazil, romildo@dme.ufcg.edu.br}\, \, and \, \  Marco A. S. Souto \thanks{M. A. S. Souto was partially supported by CNPq/Brazil
		305384/2014-7 and INCT-MAT, marco@dme.ufcg.edu.br}\,\,\,\,\,\,\vspace{2mm}
	\and {\small  Universidade Federal de Campina Grande} \\ {\small Unidade Acadêmica de Matemática} \\ {\small CEP: 58429-900, Campina Grande - PB, Brazil}\\}

\date{}

\maketitle

\begin{abstract}
	In this paper we study the existence of solution for the following class of system of elliptic equations
	$$
	\left\{
	\begin{array}{lcl}
	-\Delta u=\left(a-\int_{\Omega}K(x,y)f(u,v)dy\right)u+bv,\quad \mbox{in} \quad \Omega\\
	-\Delta v=\left(d-\int_{\Omega}\Gamma(x,y)g(u,v)dy\right)v+cu,\quad \mbox{in} \quad \Omega\\
	u=v=0,\quad \mbox{on} \quad \partial\Omega
	\end{array}
	\right.
	\eqno{(P)}
	$$
	where $\Omega\subset\R^N$ is a smooth bounded domain, $N\geq1$, and $K,\Gamma:\Omega\times\Omega\rightarrow\R$ is a nonnegative function checking some hypotheses and $a,b,c,d\in\R$. The functions $f$ and $g$ satisfy some conditions which permit to use Bifurcation Theory to prove the existence of solution for $(P)$.\vspace{0.5cm}


\end{abstract}
	
\section{Introduction and main result}

The main goal of this paper is to study the existence of positive solution for the following class of nonlocal problems 
$$
\left\{
\begin{array}{lcl}
-\Delta u=\left(a-\int_{\Omega}K(x,y)f(u,v)dy\right)u+bv,\quad \mbox{in} \quad \Omega\\
-\Delta v=\left(d-\int_{\Omega}\Gamma(x,y)g(u,v)dy\right)v+cu,\quad \mbox{in} \quad \Omega\\
u=v=0,\quad \mbox{on} \quad \partial\Omega
\end{array}
\right.
\eqno{(P)}
$$
where $\Omega\subset\R^N$ is a smooth bounded domain, $N\geq1$ and $K,\Gamma:\Omega\times\Omega\rightarrow\R$ are nonnegative functions checking some hypotheses and $a,b,c,d\in\R$. The functions $f$ and $g$ satisfy some technical conditions which will be mentioned later on.

The study of the problem $(P)$ comes from the problem to model the behavior of a species inhabiting in a smooth bounded domain $\Omega\subset\R^N$, whose the classical logistic equation is given by 
\begin{equation} \label{1}
\left\{
\begin{array}{lcl}
-\Delta u=u(\lambda-b(x)u^{p}),\quad\mbox{in} \quad\Omega\\
u=0,\quad\mbox{on}\quad\partial\Omega
\end{array}
\right.
\end{equation}
where $u(x)$ is the population density at location $x\in\Omega$, $\lambda\in\R$ is the growth rate of the species, and $b$ is a positive function denoting the carrying capacity, that is, $b(x)$ describes the limiting effect of crowding of the population.

Since (\ref{1}) is a local problem, the crowding effect of the population $u$ at $x$ only depends on the value of the population in the same point $x$. In \cite{Chipot}, for more realistic situations, Chipot has considered that the crowding effect depends also on the value of the population around of $x$, that is, the crowding effect depends on the value of integral involving the function $u$ in the ball $B_{r}(x)$ centered at $x$ of radius $r>0$. To be more precisely, in \cite{Chipot} the following nonlocal problem has been studied   
\begin{equation} \label{2}
\left\{
\begin{array}{lcl}
-\Delta u=\left(\lambda-\int_{\Omega\cap B_{r}(x)}b(y)u^{p}(y)dy\right)u,\quad\mbox{in} \quad\Omega\\
u=0,\quad\mbox{on}\quad\partial\Omega
\end{array}
\right.
\end{equation}
where $b$ is a nonnegative and nontrivial continuous function. After \cite{Chipot}, a special attention has been given for the problem 
\begin{equation} \label{3}
\left\{
\begin{array}{lcl}
-\Delta u=\left(\lambda-\int_{\Omega}K(x,y)u^{p}(y)dy\right)u,\quad\mbox{in}\quad\Omega\\
u=0,\quad\mbox{on}\quad\partial\Omega
\end{array}
\right.
\end{equation}
by supposing different conditions on $K$, see for example, Allegretto and Nistri \cite{Allegretto-Nistri}, Alves, Delgado, Souto and  Suárez \cite{Alves-Delgado-Souto-Suarez}, Chen and Shi \cite{Chen-Shi}, Corrêa, Delgado and Suárez \cite{Correa-Delgado-Suarez}, Coville \cite{Coville}, Leman, Méléard and Mirrahimi \cite{Leman-Meleard-Mirrahimi}, and Sun, Shi and Wang \cite{Sun-Shi-Wang} and their references. 

In \cite{Alves-Delgado-Souto-Suarez}, Alves, Delgado, Souto and  Suárez have considered the existence and nonexistence of solution for Problem (\ref{3}). In the paper, the authors have introduced a class $\mathcal{K}$ which is formed by functions $K:\Omega\times\Omega\rightarrow\R$ such that:

$(i)$ $K\in L^{\infty}(\Omega\times\Omega)$ and $K(x,y)\geq0$ for all $x,y\in\Omega$.

$(ii)$ If $w$ is mensurable and $\int_{\Omega\times\Omega}K(x,y)|w(y)|^{p}|w(x)|^{2}dxdy=0$, then $w=0$ a.e. in $\Omega$.

\noindent Using Bifurcation Theory and by supposing that $K$ belongs to class $\mathcal{K}$, the following result has been proved 

\begin{teo} \label{T1}
	The problem $(3)$ has a positive solution if, and only if, $\lambda>\lambda_{1}$, where $\lambda_{1}$ is the first eigenvalue of problem
	\begin{equation*}
	\left\{
	\begin{array}{lcl}
	-\Delta u=\lambda u,\quad\mbox{in}\quad\Omega\\
	u=0,\quad\mbox{on}\quad\partial\Omega.
	\end{array}
	\right.
	\end{equation*}
\end{teo} 

Motivated by \cite{Alves-Delgado-Souto-Suarez}, at least from a mathematical point of view, it seems to be interesting to ask if 
$$
\left\{
\begin{array}{lcl}
-\Delta u=\left(\lambda f(x)-\int_{\R^N}K(x,y)|u(y)|^{\gamma}dy\right)u, \quad\mbox{in}\quad\R^{N}\\
\displaystyle \lim_{|x|\to +\infty}u(x)=0,\quad u>0 \quad \text{in} \quad \R^{N}
\end{array}
\right.
\eqno{(Q)}
$$
version in $\R^N$ for (3), has a solution. This question was answered by Alves, de Lima and Souto in \cite{Alves-deLima-Souto}. In this paper, the authors study the existence of positive solution for ($Q$).

Using Bifurcation Theory and inspired by the results due to Edelson and Rumbos \cite{Edelson-Rumbos1,Edelson-Rumbos2}, Alves, de Lima and Souto, have shown that under some conditions on $K$ and $f$, problem $(Q)$ has a positive solution if, and only if, $\lambda>\lambda_{1}$, where $\lambda_{1}$  is the first eigenvalue of the linear problem
$$
	\left\{
	\begin{array}{l}
	-\Delta u=\lambda f(x)u,\quad\mbox{in}\quad\R^N\\
	\displaystyle \lim_{|x|\to +\infty}u(x)=0.
	\end{array}
	\right.
	\eqno{(AQ)}
$$
This result and your proof can be found in \cite{Alves-deLima-Souto}, where the reader can find the assumptions on $K$ and $f$. 

Motivated by \cite{Alves-Delgado-Souto-Suarez}, comes a new challenge: model the behavior of two species inhabiting in a smooth bounded domain $\Omega\subset\R^N$, analogous modeling done in the case of a species in \cite{Alves-Delgado-Souto-Suarez}. Inspired in the articles due to Corrêa and Souto \cite{Correa-Souto, Correa-Souto1} and Souto \cite{Souto}, we propose the following system to model the problem
$$
\left\{
\begin{array}{lcl}
-\Delta u=\left(a-\int_{\Omega}K(x,y)f(u,v)dy\right)u+bv,\quad\mbox{in}\quad\Omega\\
-\Delta v=\left(d-\int_{\Omega}\Gamma(x,y)g(u,v)dy\right)v+cu,\quad\mbox{in}\quad\Omega\\
u=v=0,\quad\mbox{on}\quad\partial\Omega.
\end{array}
\right.
\eqno{(P)}
$$
It is very interesting to note that in a situation where $a,b,c,d>0$, we are in a cooperative system, i.e., the two species involved mutually cooperate to their growth. If $b\cdot c<0$, we say that we are in a structure involving predator and prey. In which case $b,c<0$, there is a competition between the two species.

In the present article, well as in \cite{Alves-Delgado-Souto-Suarez}, the class $\mathcal{K}$ is formed by functions $K:\Omega\times\Omega\rightarrow\R$ such that

$i)$ $K\in L^{\infty}(\Omega\times\Omega)$ and $K(x,y)\geq0$ for all $x,y\in\Omega$.

$ii)$ If $w$ is a measurable function and $\int_{\Omega\times\Omega}K(x,y)|w(y)|^\gamma w(x)^2dxdy=0$, then $w=0$ a.e. in $\Omega$. 

The functions $K:\Omega\times\Omega\rightarrow\R$ and $\Gamma:\Omega\times\Omega\rightarrow\R$ that we are considering belong to class $\mathcal{K}$.

Related to functions $f$ and $g$, we assume that

$(f_0)$ $f,g:[0,\infty)\times[0,\infty)\rightarrow\R^+$ are continuous functions.

$(f_1)$ There exists $\epsilon>0$ such that $f(t,s)\geq\epsilon|t|^\gamma$ and $g(t,s)\geq\epsilon|s|^\gamma$, for all $t,s\in[0,\infty)$ and $\gamma>0$.

$(f_2)$ $f(pt,ps)=p^{\gamma}f(t,s)$ and $g(pt,ps)=p^{\gamma}g(t,s)$, for all $t,s\in[0,\infty)$ e $p>0$, where $\gamma>0$.

$(f_3)$ There exists $c>0$ such that $f(t,s),g(t,s)\leq c$, always that $|(t,s)|\leq 1$.

The functions $f(t,s)=|t|^{\gamma}+|s|^{\gamma-\mu}|t|^{\mu}$ and $g(t,s)=c_1|t|^{\gamma}+c_2|s|^{\gamma}$ are examples that verifies $(f_0)-(f_3)$.

The constants $a,b,c,d\in\R$, that appear in the system ($P$), forming the matrix $A=\left(\begin{array}{cc}a & b\\c & d\end{array}\right)$, which often appears in this article.

Our main results are the following:

\begin{teo} \label{TP1}
	For a matrix $A=\left(\begin{array}{cc}a & b\\c & d\end{array}\right)$ with $a,b,c,d>0$ and $\lambda>0$ its largest eigenvalue. We have that, the problem 
	$$
	\left\{
	\begin{array}{lcl}
	-\Delta u=\left(a-\int_{\Omega}K(x,y)f(u,v)dy\right)u+bv,\quad\mbox{in}\quad\Omega\\
	-\Delta v=\left(d-\int_{\Omega}\Gamma(x,y)g(u,v)dy\right)v+cu,\quad\mbox{in}\quad\Omega\\
	u,v>0,\quad\mbox{in}\quad\Omega\\
	u=v=0,\quad\mbox{on}\quad\partial\Omega.
	\end{array}
	\right.
	\eqno{(P_1)}
	$$
	has solution if, and only if, $\lambda>\lambda_1$, where $\lambda_1$ is the first eigenvalue of $(-\Delta,H_{0}^{1}(\Omega))$.
\end{teo}

In which case $f=g$ and $K=\Gamma$, we have:

\begin{teo}\label{TP2}
	Let $A=\left(\begin{array}{cc}a & b\\c & d\end{array}\right)$ a matrix such that: there is a positive and largest eigenvalue of $A$ that is the unique positive eigenvalue $\lambda$ with an eigenvector $z>0$ and $dim N(\lambda I-A)=1$. Then, the problem
	$$
	\left\{
	\begin{array}{lcl}
	-\Delta u=\left(a-\int_{\Omega}K(x,y)f(u,v)dy\right)u+bv,\quad\mbox{in}\quad\Omega\\
	-\Delta v=\left(d-\int_{\Omega}K(x,y)f(u,v)dy\right)v+cu,\quad\mbox{in}\quad\Omega\\
	u,v>0,\quad\mbox{in}\quad\Omega\\
	u=v=0,\quad\mbox{on}\quad\partial\Omega.
	\end{array}
	\right.
	\eqno{(P_2)}
	$$
	has solution for all $\lambda>\lambda_1$, where $\lambda_1$ is the first eigenvalue of $(-\Delta,H_{0}^{1}(\Omega))$.
\end{teo}

\vspace{4.3cm}

\textbf{Notations}

\begin{itemize}
	
	\item $\sigma(A)$ denotes the set of real eigenvalues of the matrix $A$.
	
	\item  $\sigma(-\Delta)$ denotes the set of eigenvalues of the operator $(-\Delta,H_{0}^{1}(\Omega))$, which has its notations and properties already well known.
	
	\item  The terms of the form $U=(u,v)$, whenever it is convenient, will be written in column matrix form
	 $U=\left(\begin{array}{c}u\\v\end{array}\right)$. Moreover, $-\Delta U=(-\Delta u,-\Delta v)$ or $-\Delta U=\left(\begin{array}{c}-\Delta u\\-\Delta v\end{array}\right)$.
 
	\item $E$ denotes the Banach space $C(\overline{\Omega})\times C(\overline{\Omega})$, with norm given by
	$$\|U\|=\|u\|_{C(\overline{\Omega})}+\|v\|_{C(\overline{\Omega})}$$
	where $U\in E$, that will always be denoted by $U=(u,v)$ or, in the column matrix form, $U=\left(\begin{array}{c}u\\v\end{array}\right)$, for $u,v\in C(\overline{\Omega})$.
	
	\item $E_1$ denotes the Banach space $C^{1}(\overline{\Omega})\times C^1(\overline{\Omega})$, with norm given by
	$$\|U\|_1=\|u\|_{C^1(\overline{\Omega})}+\|v\|_{C^1(\overline{\Omega})}$$
	where $U\in E_1$, that will always be denoted by $U=(u,v)$ or, in the column matrix form, $U=\left(\begin{array}{c}u\\v\end{array}\right)$, for $u,v\in C^1(\overline{\Omega})$.
	
	\item $z=(\alpha,\beta)>0$ or $z=\left(\begin{array}{c}\alpha\\\beta\end{array}\right)>0$, denotes that $\alpha,\beta>0$. 
	
\end{itemize}

\section{The nonlocal terms and the matricial formulation}

Suposing that $K,\Gamma\in\mathcal{K}$, as in all text, $f$ and $g$ check $(f_0)-(f_3)$, is well defined
$\phi,\psi:L^{\infty}(\Omega)\times L^{\infty}(\Omega)\rightarrow L^{\infty}(\Omega)$ given by
$$\phi_{(u,v)}(x)=\int_{\Omega}K(x,y)f(|u(y)|,|v(y)|)dy$$
and
$$\psi_{(u,v)}(x)=\int_{\Omega}\Gamma(x,y)g(|u(y)|,|v(y)|)dy.$$
Moreover, using the hypothesis of $K,\Gamma,f$ and $g$, we have the properties:

$(\phi_1)$ $t^\gamma \phi_{(u,v)}=\phi_{(tu,tv)}$ and $t^\gamma \psi_{(u,v)}=\psi_{(tu,tv)}$, for all $u,v\in L^{\infty}(\Omega)$, $t>0$ and $\gamma>0$.

$(\phi_2)$ $\|\phi_{(u,v)}\|_\infty\leq\|K\|_\infty|\Omega|\|f(|u|,|v|)\|_\infty$ and $\|\psi_{(u,v)}\|_\infty\leq\|K\|_\infty|\Omega|\|g(|u|,|v|)\|_\infty$, for all $u,v\in L^{\infty}(\Omega)$.
	
With these notations, we fix:
$$\Phi_{U}(x):=\left(\begin{array}{c}
u\phi_{(u,v)}(x)\\v\psi_{(u,v)}(x)
\end{array}\right),\quad\mbox{where}\quad U=(u,v)\in L^{\infty}(\Omega)\times L^{\infty}(\Omega).$$
	
Using the established and fixed above notations, the problem $(P_1)$ can be written in the form:
$$
\left\{
\begin{array}{lcl}
-\Delta U+\Phi_{U}(x)=AU,\quad\mbox{in}\quad\Omega\\
U>0,\quad\mbox{in}\quad\Omega\\
U=0,\quad\mbox{on}\quad\partial\Omega
\end{array}
\right.\eqno(P_3)
$$
or equivalently,
$$
\left\{
\begin{array}{lcl}
-\Delta u+\phi_{(u,v)}u=au+bv,\quad\mbox{in}\quad\Omega\\
-\Delta v+\psi_{(u,v)}u=cu+dv,\quad\mbox{in}\quad\Omega\\
u,v>0,\quad\mbox{in}\quad\Omega\\
u=v=0,\quad\mbox{on}\quad\partial\Omega.
\end{array}
\right.\eqno(P_4)
$$
Here, we recall that $U=(u,v)$ satisfies the above problem in the weak sense, if $u,v\in H_{0}^{1}(\Omega)$ and
\begin{eqnarray}
\int_{\Omega}\nabla u\nabla\varphi dx+\int_{\Omega}\phi_{(u,v)}(x)u\varphi dx=\int_{\Omega}(au+bv)\varphi dx\\
\int_{\Omega}\nabla v\nabla\eta dx+\int_{\Omega}\psi_{(u,v)}(x)v\eta dx=\int_{\Omega}(cu+dv)\eta dx
\end{eqnarray}
for all $\eta,\varphi\in H_{0}^{1}(\Omega)$.

In which case $f=g$ and $K=\Gamma$, we have $\phi_{(u,v)}=\psi_{(u,v)}$ and, consequently, $\Phi_{U}(x)=\phi(x)U$, where $\phi(x)=\phi_{(u,v)}(x)$. Thus, the problem $(P_2)$ can be written in the form
$$
\left\{
\begin{array}{lcl}
-\Delta U+\phi(x)U=AU,\quad\mbox{in}\quad\Omega\\
U>0,\quad\mbox{in}\quad\Omega\\
U=0,\quad\mbox{on}\quad\partial\Omega.
\end{array}
\right.\eqno(P_6)
$$

\section{Technical results}
	
From characteristic of our problem, it is necessary to make a technical study of matrices that include the our study.
This section is developed to present this study, which is essential in all text.	
	
\begin{lem} \label{E}
	Suppose there is a solution $U=\left(\begin{array}{c}
	u\\	v\end{array}\right)$ nontrial for the homogeneous system
	$$
	\left\{\begin{array}{cccc}
	-\Delta U &=& AU,&\quad\mbox{in}\quad\Omega\\
	U&=&0,&\quad\mbox{on}\quad\partial\Omega.
	\end{array}\right.\eqno{(Q_1)}
	$$
	Then $A$ has a real eigenvalue which is also an eigenvalue $(-\Delta,H_{0}^{1}(\Omega))$. Furthermore:
	\begin{enumerate}
		\item [i)] if $\lambda_j\in\sigma(-\Delta)\cap\sigma(A)$, for $\phi_j$ eigenfunction of $(-\Delta,H_{0}^{1}(\Omega))$ associated with the eigenvalue $\lambda_j$, we have that $z=\left(\begin{array}{c}
		\int_{\Omega}u\phi_jdx\\	\int_{\Omega}v\phi_jdx\end{array}\right)$ is eigenvector of $A$ associated with the eigenvalue $\lambda_j$.
		\item[ii)] if $\sigma(-\Delta)\cap\sigma(A)=\{\lambda_j\}$ and $dim N(A-\lambda_jI)=1$, then every solution of $(Q_1)$ is of the form $U=\phi_jz$, where $z$ is an eigenvector of $A$ associated with $\lambda_j$. Moreover, the subspace $N_{A}=\{U\in E; U\mbox{ is a solution of the problem }(Q_1)\}$ has the same dimension of the eigenspace associated with $\lambda_j$ as  eigenvalue of $(-\Delta,H_{0}^{1}(\Omega))$.
		\item[iii)] if $\sigma(-\Delta)\cap\sigma(A)=\{\lambda_j,\lambda_m\}$, $m\neq j$, then every solution of $(Q_1)$ is of the form $U=\phi_jz+\phi_mw$, where $z$ is an eigenvector of $A$ associated with $\lambda_j$ and $w$ is an eigenvector of $A$ associated with $\lambda_m$. In this case, $dim N_A$ is the sum of the dimension of the associated eigenspace with $\lambda_j$ as eigenvalue of $(-\Delta,H_{0}^{1}(\Omega))$ and the dimension of the associated eigenspace with $\lambda_m$ as eigenvalue of $(-\Delta,H_{0}^{1}(\Omega))$.
	\end{enumerate}	
\end{lem}
	
\begin{lem} \label{EE}
	Suppose that, there is a solution $U=\left(\begin{array}{c}
	u\\	v\end{array}\right)$ nonegative and nonzero for the homogeneous system
	$$
	\left\{\begin{array}{cccc}
	-\Delta U &=& AU,&\quad\mbox{in}\quad\Omega\\
	U&=&0,&\quad\mbox{on}\quad\partial\Omega.
	\end{array}\right.\eqno{(Q_1)}
	$$
	Then, $A$ has $\lambda_1$ as one of the eigenvalues, which has an eigenvector associated with positive coordinates.
\end{lem}	

\begin{cor}\label{CP}
	If $\sigma(A)=\{\mu,\lambda\}$, $\mu<\lambda$, $\lambda>0$ and $z>0$ is an eigenvector of $A$ associated to eigenvalue $\lambda$. Then, if $(Q_1)$ has $U$ as nonegative and nonzero solution, we have $\lambda=\lambda_1$ and $U=\phi_1w$, where $w$ is multiple of $z$. Moreover, we have that $U>0$ and $\frac{\partial u}{\partial\eta},\frac{\partial v}{\partial\eta}<0$ on $\partial\Omega$.
\end{cor}
	
For the reader's convenience, we present a sketch of the proofs of these above results in the appendix.

\subsection{The parameter $t$ in the homogeneous problem} \label{BI}	

Assuming that the matrix $A=\left(\begin{array}{cc}
a & b\\ c & d \end{array}\right)$ has a positive $\lambda$ eigenvalue associated to a positive $z=\left(\begin{array}{c}
\alpha\\ \beta \end{array}\right)$ eigenvector, that is, $Az=\lambda z$ with $\alpha,\beta>0$.

Our interesting is to give hypotheses about $t>0$ for that the system
$$
\left\{\begin{array}{cccc}
-\Delta U &=& tAU,&\quad\mbox{in}\quad\Omega\\
U&=&0,&\quad\mbox{on}\quad\partial\Omega.
\end{array}\right.
$$
has a space of one-dimensional solutions and a solution $U>0$ in $\Omega$.
	
Assuming the existence of a positive eigenvalue $\lambda$, consider $t=t_1=\frac{\lambda_1}{\lambda}$. And so, $t_1Az=t_1\lambda z=\lambda_1z$. Clearly, $U=\left(\begin{array}{c}
\alpha\phi_1\\ \beta\phi_1
\end{array}\right)$ is positive and satisfies $-\Delta U=t_1AU$. Therefore, the space of solutions to the problem for $t=t_1$, $N_1=N_{(t_1A)}$ has positive dimension.

In order to have a space of solutions with dimension one, we consider the situation:
\begin{itemize}
	\item if $\sigma(A)=\{\lambda,\mu\}$, with $\lambda_1\mu\neq\lambda_j\lambda$, for all $j>1$. In this case, $\sigma(t_1A)=\{t_1\mu,\lambda_1\}$. As $t_1\mu\neq\lambda_j$, for all $j>1$, follow that $dim N_1=1$. 
\end{itemize}

It is easy to see that if $\lambda>\mu$, the above condition is always satisfied: $t_1\mu<t_1\lambda=\lambda_1<\lambda_j$, for all $j>1$.  

The situation here descript is utilized in the Theorems \ref{TP1} and \ref{TP2}.
	
On the other hand, since we will make use of the global bifurcation theorem, note that if $A$ has two positive eigenvalues $\lambda$ and $\mu$, and each is associated positive eigenvectors $z$ and $w$, respectively, then $dim N_{(tA)}=1$, for $t=t_1$ and $t=s_1$, where $t_1=\lambda_1/\lambda$ and $s_1=\lambda_1/\mu$. Moreover, $0<z\phi_1\in N_{(t_1A)}$ and $0<w\phi_1\in N_{(s_1A)}$. Thus, a bifurcation can starts at $t=t_1$ and may finish in $t=s_1$. We must avoid this situation.
	
Therefore, the hypothesis on $A$ is that this matrix has at least one positive eigenvalue $\lambda$ with $dim N(A-\lambda I)=1$, and it is associated with a positive eigenvector $z$. Moreover, if $A$ has another positive eigenvalue $\mu$, must be associated with an eigenvector $w=\left(\begin{array}{c}
\alpha_2\\ \beta_2\end{array}\right)$ with $\alpha_2\beta_2<0$.

By Lemma \ref{sinal}, a bifurcation with positive solutions should start in $t=t_1$.

\begin{ob}\textbf{Recalling linear algebra of the matrices $2\times 2$}

For a matrix $A=\left(\begin{array}{cc}a & b\\c & d\end{array}\right)$ with $a,b,c,d>0$, it is possible to prove that $\sigma(A)=\{\mu,\lambda\}$ and $\lambda>\mu$ with $\lambda>0$. It is well know that there exist $\left(\begin{array}{c}\alpha\\\beta\end{array}\right)$ eigenvector of $A$ associated with the eigenvalue $\lambda$ with $\alpha,\beta>0$ and $\left(\begin{array}{c}\alpha_1\\\beta_1\end{array}\right)$ eigenvector of $A$ associated with the eigenvalue $\mu$ with $\alpha_1\beta_1<0$.	
\end{ob} 

\section{Comments on the solutions operators}
	
We intend to prove the existence of positive solution for $(P_1)$ and $(P_2)$ by using the classical bifurcation result due to Rabinowitz, see \cite{Rabinowitz}. To this end, we recall that there exists $c_\infty=c_\infty(\Omega)>0$ such that: for each $h\in L^{\infty}(\Omega)$, there is only $\omega\in C^{1}(\overline{\Omega})$ satisfying:
$$
\left\{\begin{array}{lcl}
-\Delta\omega=h(x),\mbox{ in }\Omega\\
\omega=0,\mbox{ on }\partial\Omega
\end{array}
\right.
$$
and
$$\|\omega\|_{C^{1}(\overline{\Omega})}\leq c_\infty\|h\|_\infty.$$
We use the property up freely in the construction of elementary properties of operators who build below.

Considering the solution operator $S:E\rightarrow E_1$, given by
$$S(u,v)=(u_1,v_1)\Leftrightarrow\left\{\begin{array}{c}
-\Delta u_1=au+bv,\quad\mbox{in}\quad\Omega\\
-\Delta v_1=cu+dv,\quad\mbox{in}\quad\Omega\\
u_1=v_1=0,\quad\mbox{on}\quad\partial\Omega
\end{array}\right.$$
or, equivalently, in the matricial form
$$S(U)=U_1\Leftrightarrow\left\{\begin{array}{c}
-\Delta U_1=AU,\quad\mbox{in}\quad\Omega\\U_1=0,\quad\mbox{on}\quad\partial\Omega
\end{array}\right.$$
where $U=\left(\begin{array}{c}u\\v\end{array}\right)$ and $U_1=\left(\begin{array}{c}u_1\\v_1\end{array}\right)$. We have that, $S$ is well defined, is linear and verifies
$$\|S(U)\|_1\leq C\|U\|,\mbox{ for all }U\in E.$$
Moreover, using the Schauder embedding $S:E\rightarrow E$ is a compact operator.

On the other hand, setting the nonlinear operator $G:E\rightarrow E_1$ given by
$$G(u,v)=(u_1,v_1)\Leftrightarrow\left\{\begin{array}{c}
-\Delta u_1+\phi_{(u,v)}(x)u=0,\quad\mbox{in}\quad\Omega\\
-\Delta v_1+\psi_{(u,v)}(x)v=0,\quad\mbox{in}\quad\Omega\\
u_1=v_1=0,\quad\mbox{on}\quad\partial\Omega
\end{array}\right.$$
or, equivalently, in matricial form
$$G(U)=U_1\Leftrightarrow\left\{\begin{array}{c}
-\Delta U_1+\Phi_U(x)=0,\quad\mbox{in}\quad\Omega\\U_1=0,\quad\mbox{on}\quad\partial\Omega
\end{array}\right.$$
where $\Phi_{U}(x)=\left(\begin{array}{c}
u\phi_{(u,v)}(x)\\v\psi_{(u,v)}(x)
\end{array}\right)$ and $U=(u,v)$. We have, clearly, that $G$ is well defined, it is continuous and checks
$$\|G(U)\|_1\leq C(\|\phi_{(u,v)}\|_\infty+\|\psi_{(u,v)}\|_\infty)\|U\|,\mbox{ for all }U\in E.$$
Using again the Schauder embedding, we have that $G:E\rightarrow E$ is compact. Moreover, from $(f_2)-(f_3)$ and $(\phi_2)$, is possible to verify that
$$G(U)=o(\|U\|).$$
 
\section{Proof of Theorem \ref{TP1}}

In order to prove Theorem \ref{TP1} via bifurcation theory, it is necessary to introduce a parameter $t>0$ in the problem ($P_1$) and prove the lemma below:

\begin{lem} \label{LP}
	For a matrix $A=\left(\begin{array}{cc}a & b\\c & d\end{array}\right)$ with $a,b,c,d>0$ and $\lambda>0$ its largest eigenvalue. We have that, the problem 
	$$
	\left\{
	\begin{array}{lcl}
	-\Delta U+\Phi_{U}(x)=tAU,\quad\mbox{in}\quad\Omega\\
	U>0,\quad\mbox{in}\quad\Omega\\
	U=0,\quad\mbox{on}\quad\partial\Omega
	\end{array}
	\right.\eqno(P_7)
	$$
	has solution if, and only if, $t>t_1$, where  $t_1=\frac{\lambda_1}{\lambda}$ and  $\lambda_1$ is the first eigenvalue of $(-\Delta,H_{0}^{1}(\Omega))$.
\end{lem}

To prove the lemma above, it is necessary to note that: using the definitions of $S$ and $G$, it is easy to check that $(t,U)\in\R\times E$ solves $(P_7)$ if, and
only if,
$$U=F(t,U)=tS(U)+G(U).$$

In the sequel, we will apply the following result due to Rabinowitz \cite{Rabinowitz}, to prove the Lemma \ref{LP}.

\begin{teo}\textbf{(Global bifurcation)}\label{bifur}
	Let $E$ be a Banach space. Suppose that $S$ is a compact linear operator
	and $t^{-1}\in \sigma(S)$ has odd algebraic multiplicity . If $G$ is a compact operator and
	$$
	\lim_{\|u\|\to 0}\frac{G(u)}{\|u\|}=0,
	$$
	then  the set
	$$
	\Sigma=\overline{\{(t,u)\in\R\times E:u=tS(u)+G(u),u\neq0\}}
	$$
	has a closed connected component $\mathcal{C}=\mathcal{C}_{t}$ such that $(t,0)\in\mathcal{C}$ and
	
	(i) $\mathcal{C}$ is unbounded in $\R\times E$, 
	or
	
	(ii) there exists $\hat{t}\neq t$, such that $(\hat{t},0)\in\mathcal{C}$ and $\hat{t}^{-1}\in\sigma(S)$.
\end{teo}

By study done in Subsection \ref{BI}, an eigenfunction $U_1$ associated with eigenvalue $t_1=\lambda_1/\lambda$ of the linear problem can be chosen positive. In addition, $t_1^{-1}$ is an eigenvalue of multiplicity 1 for $S$. From global bifurcation theorem, there exists a closed connected component $\mathcal{C}=\mathcal{C}_{t_1}$ of solutions for $(P_7)$, which satisfies $(i)$ or $(ii)$. We claim that $(ii)$ does not occur. In order to show this claim, we need the two lemmas below

\begin{lem}\label{sinal}
	There exists $\delta>0$ such that, if $(t,U)\in\mathcal{C}$ with $|t-t_1|+\|U\|<\delta$ and $U\neq0$, then $U$ has defined signal, that is,
	$$
	U(x)>0, \quad \forall x\in\Omega\quad \mbox{or} \quad U(x)<0,\quad \forall x\in\Omega.
	$$
\end{lem}

\begin{dem}
	It is enough to prove that for any two sequences $(U_n) \subset E$ and $t_{n}\to t_1$ with 
	$$
	U_n\neq0,\quad \|U_n\|\to 0 \quad \mbox{and } \quad U_{n}=F(t_n,U_n)=t_{n}S(U_n)+G(U_n),
	$$
	$U_n$ has defined signal for $n$ large enough. 
	
	Setting $W_n=U_n/\|U_n\|$, we have that
	$$W_n=t_nS(W_n)+\frac{G(U_n)}{\|U_n\|}=t_nS(W_n)+o_n(1).$$
	From compactness of the operator $S$, we can assume that $(S(W_n))$ is convergent. Then, $W_n\to W$ in $E$ for some $W \in E$ with $\|W\|=1$. Consequently, 
	$$
	\left\{\begin{array}{cccc}
	-\Delta W &=& t_1AW,&\mbox{in}\quad\Omega\\
	W&=&0,&\mbox{on}\quad\partial\Omega.
	\end{array}\right.
	$$	
	Once that $W\neq0$, we have, by Lemmas \ref{E} and \ref{EE}, that
	$$W(x)>0\mbox{ or }W(x)<0,\mbox{ for all }x\in\Omega.$$
	Therefore, without loss of generality, $W>0$ in $\Omega$, and consequently $W_n>0$ in $\Omega$ for $n$ large enough. Once  $U_n$ and $W_n$ has the same signal, we have that $U_n$ is also positive, this completes to proof.
\end{dem}

It is easy to check that if $(t,U)\in\Sigma$, the pair $(t,-U)\in\Sigma$. From maximum principle arguments used in \cite{Alves-Delgado-Souto-Suarez} and positivity of $a,b,c$ and $d$, we can decompose  $\mathcal{C}$ in $\mathcal{C}^+\cup\mathcal{C}^-$, where
$$\mathcal{C}^+=\{(t,U)\in\mathcal{C};U>0\}\cup\{(t_1,0)\}$$
and
$$\mathcal{C}^-=\{(t,U)\in\mathcal{C};U<0\}\cup\{(t_1,0)\}.$$
Observed that, $\mathcal{C}^-=\{(t,U)\in\mathcal{C};(t,-U)\in\mathcal{C}^+\}$, $\mathcal{C}^+\cap\mathcal{C}^-=\{(t_1,0)\}$ and $\mathcal{C}^+$ is
unbounded if, and only if, $\mathcal{C}^-$ is unbounded.

Now, we are able to prove that $(ii)$ does not hold.
\begin{lem} \label{CUB}
	$\mathcal{C}^{+}$ is unbounded.
\end{lem} 

\begin{dem}
	Suppose by contradiction that $\mathcal{C}^{+}$ is bounded. Then, $\mathcal{C}$ is also bounded. From global bifurcation theorem, there exists $(\hat{t},0)\in\mathcal{C}$ , where $\hat{t}\neq t_1$ e $\hat{t}^{-1}\in\sigma(S)$.
	
	Hence, without loss of generality, there exist $(t_n,U_n)\subset\mathcal{C}^+$ with $t_n\rightarrow\hat{t}$ such that
	$$U_n\neq0,\quad\|U_n\|\rightarrow0\mbox{ and }U_n=F(t_n,U_n).$$
	
	Setting $W_n=U_n/\|U_n\|$, similar to what was done in the previous lemma, there exists $W\in E$ with $W_n\rightarrow W$ in $E$, where $W\neq0$, $W\geq0$ and satisfies
	$$
	\left\{\begin{array}{cccc}
	-\Delta W &=& (\hat{t}A)W,&\mbox{in}\quad\Omega\\
	W&=&0,&\mbox{in}\quad\partial\Omega.
	\end{array}\right.
	$$
	From Corollary \ref{CP}, $\hat{t}\lambda=\lambda_1$ and, consequently, $\hat{t}=t_1$, which is impossible. This proves the lemma.
\end{dem}
	
From previous lemma, the connected component $\mathcal{C}^{+}$ is unbounded. Now, our goal is to show that this component
intersects any hyperplane $\{t\}\times E$, for $t>t_{1}$. To see this, we need of the following a priori estimate

\begin{lem} \label{EP}\textbf{(A priori estimate)} For any $\Lambda>0$, there exists $R>0$ such that, if $(t,U)\in\mathcal{C}^{+}$ and $t \in [0,\Lambda]$, then $\|U\|\leq R$. 
\end{lem} 

\begin{dem}
	Setting by $\|\cdot\|_H$, the norm in $H=H_{0}^{1}(\Omega)\times H_{0}^{1}(\Omega)$, given by
	$$\|U\|_H=\|u\|_{H_{0}^{1}(\Omega)}+\|v\|_{H_{0}^{1}(\Omega)}$$
	where $U=\left(\begin{array}{c}u\\v\end{array}\right)$ with $u,v\in H_{0}^{1}(\Omega)$.
	
	We start showing an a priori estimate on the $H$ space:
	
	\begin{Af}
		Given $\Lambda>0$, there exists $R>0$ such that: if $(t,U)\in\mathcal{C}^+$ and $t\leq\Lambda$, then $\|U\|_H\leq R$.
	\end{Af}
	Indeed, if the claim does not hold, there are $(U_n)\subset H$ and $(t_n)\subset[0,\Lambda]$ such that,
	$$\|U_n\|_H\rightarrow\infty\quad\mbox{and}\quad U_n=F(t_n,U_n).$$
	Consider $W_n=U_n/\|U_n\|_H$, where $W_n=(\overline{u}_n,\overline{v}_n)$ for $\overline{u}_n=u_n/\|U_n\|_H$ and $\overline{v}_n=v_n/\|U_n\|_H$. Thus,
	$$
	\int_{\Omega}\nabla\overline{u}_n\nabla\varphi dx+\int_{\Omega}\phi_{(u_n,v_n)}(x)\overline{u}_n\varphi dx=t_n\int_{\Omega}(a\overline{u}_n+b\overline{v}_n)\varphi dx
	$$
	and
	$$
	\int_{\Omega}\nabla\overline{v}_n\nabla\eta dx+\int_{\Omega}\psi_{(u_n,v_n)}(x)\overline{v}_n\eta dx=t_n\int_{\Omega}(c\overline{u}_n+d\overline{v}_n)\eta dx
	$$
	for all $\varphi,\eta\in H_{0}^{1}(\Omega)$. Once $(W_n)$ is bounded in $H$, without loss of generality, we can suppose that there is $W\in H$, $W=(u,v)$, such that
	\begin{eqnarray}
	\overline{u}_n\rightharpoonup u\mbox{ in }H_{0}^{1}(\Omega),\overline{u}_n\rightarrow u\mbox{ in }L^2(\Omega)\mbox{ and }\overline{u}_n(x)\rightarrow u(x)\mbox{ a.e. in }\Omega\\
	\overline{v}_n\rightharpoonup v\mbox{ in }H_{0}^{1}(\Omega),\overline{v}_n\rightarrow v\mbox{ in }L^2(\Omega)\mbox{ and }\overline{v}_n(x)\rightarrow v(x)\mbox{ a.e. in }\Omega.
	\end{eqnarray}
	
	For $\varphi=\displaystyle{\frac{\overline{u}_n}{\|U\|_H^{\gamma}}}$ and $\eta=\displaystyle{\frac{\overline{v}_n}{\|U\|_H^{\gamma}}}$ as functions of test, and recaling that $t^\gamma\phi_{(u_n,v_n)}=\phi_{(tu_n,tv_n)}$ and also $t^\gamma\psi_{(u_n,v_n)}=\psi_{(tu_n,tv_n)}$, for all $t>0$, getting
	\begin{eqnarray}
	\frac{1}{\|U_n\|_{H}^{\gamma}}\|\overline{u}_n\|^{2}_{H_{0}^{1}(\Omega)}+\int_{\Omega}\phi_{(\overline{u}_n,\overline{v}_n)}\overline{u}_n^2dx=t_n\int_{\Omega}(a\overline{u}_n+b\overline{v}_n)\frac{\overline{u}_n}{\|U_n\|_{H}^{\gamma}}dx\\
	\frac{1}{\|U_n\|_{H}^{\gamma}}\|\overline{v}_n\|^{2}_{H_{0}^{1}(\Omega)}+\int_{\Omega}\psi_{(\overline{u}_n,\overline{v}_n)}\overline{v}_n^2dx=t_n\int_{\Omega}(c\overline{u}_n+d\overline{v}_n)\frac{\overline{v}_n}{\|U_n\|_{H}^{\gamma}}dx.
	\end{eqnarray}
	Therefore, using Hölder inequality and $(6)-(7)$ in $(8)-(9)$, 
	\begin{equation}
	\lim_{n\rightarrow\infty}\int_{\Omega}\phi_{(\overline{u}_n,\overline{v}_n)}\overline{u}_n^2dx=\lim_{n\rightarrow\infty}\int_{\Omega}\psi_{(\overline{u}_n,\overline{v}_n)}\overline{v}_n^2dx=0.
	\end{equation}
	From Fatou lemma,
	\begin{eqnarray}
	\int_{\Omega}\phi_{(u,v)}u^2dx\leq\lim_{n\rightarrow\infty}\int_{\Omega}\phi_{(\overline{u}_n,\overline{v}_n)}\overline{u}_n^2dx=0\\
	\int_{\Omega}\psi_{(u,v)}v^2dx\leq\lim_{n\rightarrow\infty}\int_{\Omega}\psi_{(\overline{u},\overline{v}_n)}\overline{v}_n^2dx=0.
	\end{eqnarray}
	And so,
	\begin{equation}
	\int_{\Omega\times\Omega}K(x,y)f(|u(y)|,|v(y)|)|u(x)|^2dxdy=\int_{\Omega\times\Omega}\Gamma(x,y)g(|u(y)|,|v(y)|)|v(x)|^2dxdy=0.
	\end{equation}
	Thus, by $(f_1)$,
	\begin{eqnarray}
	0\leq\epsilon\int_{\Omega\times\Omega}K(x,y)|u(y)|^{\gamma}|u(x)|^2dxdy\leq\int_{\Omega\times\Omega}K(x,y)f(|u(y)|,|v(y)|)|u(x)|^2dxdy=0\\
	0\leq\epsilon\int_{\Omega\times\Omega}\Gamma(x,y)|v(y)|^{\gamma}|v(x)|^2dxdy\leq\int_{\Omega\times\Omega}\Gamma(x,y)g(|u(y)|,|v(y)|)|v(x)|^2dxdy=0.
	\end{eqnarray}
	and, consequently
	\begin{equation}
	\int_{\Omega\times\Omega}K(x,y)|u(y)|^{\gamma}|u(x)|^2dxdy=\int_{\Omega\times\Omega}\Gamma(x,y)|v(y)|^{\gamma}|v(x)|^2dxdy=0.
	\end{equation}
	Since that $K,\Gamma\in\mathcal{K}$, we have that $u=v=0$. Hence, $(\overline{u}_n)$ and $(\overline{v}_n)$ converge to 0 in $L^{2}(\Omega)$. Considering $\varphi=\overline{u}_n$ and $\eta=\overline{v}_n$ as test function, we get that
	\begin{eqnarray}
	\int_{\Omega}|\nabla\overline{u}_n|^2dx+\int_{\Omega}\phi_{(u_n,v_n)}\overline{u}_n^2dx=t_n\int_{\Omega}(a\overline{u}_n+b\overline{v}_n)\overline{u}_ndx\\
	\int_{\Omega}|\nabla\overline{v}_n|^2dx+\int_{\Omega}\psi_{(u_n,v_n)}\overline{v}_n^2dx=t_n\int_{\Omega}(c\overline{u}_n+d\overline{v}_n)\overline{v}_ndx.
	\end{eqnarray}
	As $(t_n)$ is bounded by $\Lambda$,
	\begin{eqnarray}
	\int_{\Omega}|\nabla\overline{u}_n|^2dx\leq \Lambda\left[a\int_{\Omega}|\overline{u}_n|^2dx+b\int_{\Omega}|\overline{u}_n\overline{v}_n|dx\right]\\
	\int_{\Omega}|\nabla\overline{v}_n|^2dx\leq \Lambda\left[c\int_{\Omega}|\overline{u}_n\overline{v}_n|dx+d\int_{\Omega}|\overline{v}_n|^2dx\right].
	\end{eqnarray}
	Consequently, $\|W_n\|_H\rightarrow0$. This contradictis the fact that $\|W_n\|_H=1$ for all $n\in\N$, and the lemma follows.
	
	Since $(U_n)$ is bounded in $H$, iteration arguments imply that $(U_n)$ is bounded in $L^{\infty}(\Omega)\times L^{\infty}(\Omega)$, and the proof is done.
\end{dem}

\noindent {\bf Conclusion of the proof of Lemma \ref{LP} and proof of the Theorem \ref{TP1}}

From Lemma \ref{EP}, for all $t>t_1$, we have that $(\{t\}\times E)\cap\mathcal{C}^+\neq\emptyset$, that is, $\mathcal{C}^+$ crosses the hyperplane $\{t\}\times E$. Indeed, otherwise there $\Lambda>t_1$ such that $\mathcal{C}^+$ does not cross the hyperplane $\{\Lambda\}\times E$, thus by Lemma \ref{EP} there exists $R>0$ such that $(t,U)\in\mathcal{C}^+$, $t\in[0,\Lambda]$, and $\|U\|\leq R$. Therefore, $\mathcal{C}^+$ would be bounded, which contradicts the Lemma \ref{CUB}. 

To finalize the proof of Lemma \ref{LP}, we must show that there is no solution for $(P_7)$ when $t \leq t_{1}=\frac{\lambda_1}{\lambda}$. Indeed, arguing by contradiction, if $(t,U)$ is a solution of $(P_7)$, with $t\leq t_1$ and $U=(u,v)>0$, we have
$$
\left\{
\begin{array}{lcl}
-\Delta u+\phi_{(u,v)}u=t[au+\frac{b}{\sigma}(\sigma v)],\quad\mbox{in}\quad\Omega\\
-\Delta (\sigma v)+\psi_{(u,v)}(\sigma v)=t[(c\sigma)u+d(\sigma v)],\quad\mbox{in}\quad\Omega\\
u=v=0,\quad\mbox{on}\quad\partial\Omega
\end{array}
\right.
$$
for all $\sigma>0$. In particular, if $\sigma^2=\frac{b}{c}$, may we fix $w=\sigma v$ and observe that $\hat{b}:=\frac{b}{\sigma}=c\sigma$. And so, $\hat{U}=(u,w)\in E$ is solution of the problem 
$$
\left\{
\begin{array}{lcl}
-\Delta u+\phi_{(u,v)}u=t[au+\hat{b}w],\quad\mbox{in}\quad\Omega\\
-\Delta w+\psi_{(u,v)}w=t[\hat{b}u+dw],\quad\mbox{in}\quad\Omega\\
u,w>0,\quad\mbox{in}\quad \Omega\\
u=w=0,\quad\mbox{on}\quad\partial\Omega
\end{array}
\right.
$$
Since $A_0=\left(\begin{array}{cc}a & \hat{b}\\ \hat{b} & d\end{array}\right)$ is a simmetric matrix, we know that
\begin{equation}\label{propriedade}
\mu|z|^2\leq\left<A_0z,z\right>\leq \lambda|z|^2,\quad \mbox{for all}\quad z\in\R^2.
\end{equation}
On the other hand,
\begin{equation}
\int_{\Omega}\left<tA_0\left(\begin{array}{c}u\\w\end{array}\right),\left(\begin{array}{c}u\\w\end{array}\right)\right>dx=\int_{\Omega}|\nabla u|^2+\phi_{(u,v)}(x)u^2dx+\int_{\Omega}|\nabla w|^2+\psi_{(u,v)}(x)w^2dx>\int_{\Omega}(|\nabla u|^2+|\nabla w|^2)dx
\end{equation}
and, consequently by (\ref{propriedade}),
\begin{equation}
\int_{\Omega}(|\nabla u|^2+|\nabla w|^2)dx<t\lambda\int_{\Omega}(|u|^2+|w|^2)dx.
\end{equation}
On the other hand, by Poincaré's inequality
\begin{equation}
\lambda_1\int_{\Omega}(|u|^2+|w|^2)dx<t\lambda\int_{\Omega}(|u|^2+|w|^2)dx.
\end{equation}
Hence, $t>\frac{\lambda_1}{\lambda}$ which is a contradiction. This proves the lemma.	

In relation to Theorem \ref{TP1}, by Lemma \ref{LP}, it is clear that ($P_1$) has solution if, and only if, $t_1=\frac{\lambda_1}{\lambda}<1$. Therefore, ($P_1$) has solution if, and only if, $\lambda>\lambda_1$. This proves the theorem.

\section{Proof of Theorem \ref{TP2}}

As was done in the Theorem \ref{TP1}, to prove the Theorem \ref{TP2} via bifurcation theory, it is necessary to introduce a parameter $t>0$ in the problem ($P_2$) and prove the lemma below:

\begin{lem}\label{LP2}
	Let $A=\left(\begin{array}{cc}a & b\\c & d\end{array}\right)$ a matrix such that: there is a positive and largest eigenvalue of $A$ that is the unique positive eigenvalue $\lambda$ with an eigenvector $z>0$ and $dim N(\lambda I-A)=1$. Then, the problem
	$$
	\left\{
	\begin{array}{lcl}
	-\Delta U+\phi(x)U=tAU,\quad\mbox{in}\quad\Omega\\
	U>0,\quad\mbox{in}\quad\Omega\\
	U=0,\quad\mbox{on}\quad\partial\Omega
	\end{array}
	\right.\eqno(P_8)
	$$
	has solution for all $t>t_1$, where  $t_1=\frac{\lambda_1}{\lambda}$ and  $\lambda_1$ is the first eigenvalue of $(-\Delta,H_{0}^{1}(\Omega))$.
\end{lem}

To prove the lemma above, it is necessary to note that: using the definitions of $S$ and $G$, it is easy to check that $(t,U)\in\R\times E$ solves $(P_8)$ if, and
only if,
$$U=F(t,U)=tS(U)+G(U).$$

In the sequel, we will apply again the result due to Rabinowitz \cite{Rabinowitz}, for prove the Lemma \ref{LP2}.
	
Again, by study done in Subsection \ref{BI}, a eigenfunction $U_1$ associated with eigenvalue $t_1=\lambda_1/\lambda$ of the linear problem can be chosen positive. In addition, $t_1^{-1}$ is an eigenvalue of multiplicity 1 for $S$. From global bifurcation theorem (Theorem \ref{bifur}), there exists a closed connected component $\mathcal{C}=\mathcal{C}_{t_1}$ of solutions for $(P_8)$, which satisfies $(i)$ or $(ii)$, of the Theorem \ref{bifur}. We claim that $(ii)$ does not occur. In order to show this claim, we need the lemma below
\begin{lem}\label{sinal2}
	There exists $\delta>0$ such that, if $(t,U)\in\mathcal{C}$ with $|t-t_1|+\|U\|<\delta$ and $U\neq0$, then $U$ has defined signal, that is,
	$$
	U(x)>0, \quad \forall x\in\Omega\quad \mbox{or} \quad U(x)<0,\quad \forall x\in\Omega.
	$$
\end{lem}
that is, the same Lemma \ref{sinal}, but now on the problem $(P_8)$. The demonstration is absolutely analogous with trivial modifications.

A trivial fact in the previous case is the decomposition of $\mathcal{C}$ in $\mathcal{C}^+\cup\mathcal{C}^-$. Here, to see such decomposition, a special attention is required.

In order to achieve this decomposition, it is necessary to introduce an auxiliary operator, whose properties are similar to the Laplacian operator.

\textbf{Auxiliary operator}
	
By fixing $\psi\in L^{\infty}(\overline{\Omega})$, the solution operator $S_L:L^{2}(\Omega)\rightarrow L^{2}(\Omega)$ such that $S_L(v)=u$, where $u$ is the unique weak solution for the linear problem
\begin{equation}
\left\{
\begin{array}{lcl}
L(u)=v,\quad\mbox{in}\quad\Omega\\
u=0,\quad\mbox{on}\quad\partial\Omega
\end{array}
\right.
\end{equation} 
where $L(u)=-\Delta u+\psi(x)u$. This solution operator is compact self-adjoint, then by spectral theory there exists a complete orthonormal basis $\{\phi_n\}$ of $L^{2}(\Omega)$ and a corresponding sequence of positive real numbers $\{\lambda_{n}\}$ with $\lambda_{n}\rightarrow\infty$ as $n\rightarrow\infty$ such that
$$
0<\lambda_{1}\leq\lambda_{2}\leq...\leq \lambda_n \leq .....
$$
and
$$
\left\{
\begin{array}{lcl}
L(\phi_n)=\lambda_n \phi_n,\quad\mbox{in}\quad\Omega\\
\phi_n=0,\quad\mbox{on}\quad\partial\Omega.
\end{array}
\right.
$$
Moreover, using Lagrange multiplier it is possible to prove the following characterization for $\lambda_1$ 
\begin{equation}
\lambda_1=\inf_{v\in H^1(\Omega)\setminus\{0\}}\frac{\int_{\Omega}[|\nabla v|^2+\psi(x)v^2]dx}{\int_{\Omega}v^2dx}.
\end{equation}
The above identity is crucial to show that $\lambda_{1}$ is a simple eigenvalue and that a corresponding eigenfunction $\phi_{1}$ can be chosen positive in $\Omega$. Note that the Lemma \ref{E} is valid, replacing $-\Delta U$ by $LU$ and $\sigma(-\Delta)$ by $\sigma(L)$, where $LU=(L(u),L(v))$.

With the notations and properties introduced above, we have the following version of the Hopf's Lemma in matricial format:

\begin{lem}\label{operadorauxiliar}
	If the problem
	$$
	\left\{
	\begin{array}{lcl}
	LU=AU,\quad\mbox{in}\quad\Omega\\
	U=0,\quad\mbox{on}\quad\partial\Omega
	\end{array}
	\right.
	$$
	has solution $U$ with $U\geq0$ and $U\neq0$, then $\sigma(L)\cap\sigma(A)\neq\emptyset$. Moreover, $U>0$ in $\Omega$ and $\frac{\partial u}{\partial\eta},\frac{\partial v}{\partial\eta}<0$ on $\partial\Omega$. 
\end{lem}

\begin{dem}
	Using the same arguments of the Lemma \ref{E}, we conclude that $\sigma(L)\cap\sigma(A)\neq\emptyset$. Furthermore, as there is a positive and largest eigenvalue of $A$ that is the unique positive eigenvalue $\lambda$ with an eigenvector $z>0$ and $dim N(\lambda I-A)=1$, we have, by analogous argument to the Corollary \ref{CP}, $\lambda=\lambda_1$ and $U=\phi_1w$,  where $w$ is multiple of $z$. Moreover, we have that $U>0$ and $\frac{\partial u}{\partial\eta},\frac{\partial v}{\partial\eta}<0$ on $\partial\Omega$. This completes the proof.
\end{dem}

We can now prove the theorem:

\begin{lem} Consider the sets 
	$$
	\mathcal{C}^{+}=\{(t,U)\in\mathcal{C}:U(x)>0,\quad \forall x\in\Omega \}\cup\{(t_{1},0)\}
	$$
	and
	$$
	\mathcal{C}^{-}=\{(t,U)\in\mathcal{C}:U(x)<0,\quad \forall x\in\Omega\}\cup\{(t_{1},0)\}.
	$$
	Then, 
	\begin{equation} \label{Componente}
	\mathcal{C}=\mathcal{C}^{+}\cup\mathcal{C}^{-}.
	\end{equation}
	Moreover, note that   
	$\mathcal{C}^{-}=\{(t,U)\in\mathcal{C}:(t,-U)\in\mathcal{C}^{+}\}$, $\mathcal{C}^{+}\cap\mathcal{C}^{-}=\{(t_{1},0)\}$ and $\mathcal{C}^{+}$ is unbounded if, and only if, $\mathcal{C}^{-}$ is also unbounded.
\end{lem}

\begin{dem}
	Of course, the proof is complete, showing that $\mathcal{C}^{+}$ is closed and open. For $(t,U)\in\overline{\mathcal{C}^{+}}$, we have $U\neq0$ and $U\geq0$, where $U=(u,v)$. As,
	\begin{equation}
	\left\{
	\begin{array}{lcl}
	-\Delta u+\phi_{(u,v)}u=au+bv,\quad\mbox{in}\quad\Omega\\
	-\Delta v+\phi_{(u,v)}v=cu+dv,\quad\mbox{in}\quad\Omega\\
	u=v=0,\quad\mbox{on}\quad\partial\Omega
	\end{array}
	\right.
	\end{equation}
	we get,
	\begin{equation}
	\left\{
	\begin{array}{lcl}
	LU=tAU,\quad\mbox{in}\quad\Omega\\
	U=0,\quad\mbox{on}\quad\partial\Omega
	\end{array}
	\right.
	\end{equation}
	for $L(w):=-\Delta w+\phi(x)w$, where $\phi(x)=\phi_{(u,v)}(x)$. 
	Therefore, by Lemma \ref{operadorauxiliar}, we obtain $U>0$ in $\Omega$ and, consequently, $\mathcal{C}^{+}$ is closed. Now, for $(t,U)\in\mathcal{C}^{+}$, we have, by Lemma \ref{operadorauxiliar}, that $U>0$ and $\frac{\partial u}{\partial\eta},\frac{\partial v}{\partial\eta}<0$ on $\partial\Omega$. Therefore, by Hopf's Lemma, $(t,U)\in int\mathcal{C}^{+}$. This proves the lemma.
\end{dem}

Following the same steps of the previous section, we should prove:
\begin{lem}\label{CUB2}
	$\mathcal{C}^{+}$ is unbounded.
\end{lem}
But, the proof is absolutely analogous, with some trivial modifications, to the proof of the Lemma \ref{CUB}, so not we shall prove. The same can be said of the a priori estimate:

\begin{lem}\label{EP2}\textbf{(A priori estimate)} For any $\Lambda>0$, there exists $R>0$ such that, if $(t,U)\in\mathcal{C}^{+}$ and $t \in [0,\Lambda]$, then $\|U\|\leq R$. 
\end{lem} 

\noindent {\bf Conclusion of the proof of Lemma \ref{LP2} and proof of the Theorem \ref{TP2}}

From Lemma \ref{EP2}, for all $t>t_1$, we have that $(\{t\}\times E)\cap\mathcal{C}^+\neq\emptyset$, that is, $\mathcal{C}^+$ crosses the hyperplane $\{t\}\times E$. Indeed, otherwise there $\Lambda>t_1$ such that $\mathcal{C}^+$ does not cross the hyperplane $\{\Lambda\}\times E$, thus by Lemma \ref{EP2} there exists $R>0$ such that $(t,U)\in\mathcal{C}^+$, $t\in[0,\Lambda]$ and $\|U\|\leq R$. Therefore, $\mathcal{C}^+$ would be bounded, which contradicts the Lemma \ref{CUB2}. 

In relation to Theorem \ref{TP2}, by Lemma \ref{LP2}, it is clear that ($P_2$) has solution if $t_1=\frac{\lambda_1}{\lambda}<1$. Therefore, ($P_2$) has solution if $\lambda>\lambda_1$. This proves the theorem.

\vspace{1.3cm}

\section{Appendix}

This section is dedicated to present some details that could be removed from the text without prejudice to the understanding of the content.

\subsection{Proof of the Lemmas \ref{E} and \ref{EE}}

\begin{dem}(\textbf{Lemma \ref{E}})
	
	Suppose that, $\int_{\Omega}u\phi_jdx\neq0$. Multiplying the equations in $(Q_1)$ by $\phi_j$ and integrating on $\Omega$, we get:
	$$\lambda_j\int_{\Omega}u\phi_jdx=\int_{\Omega}\nabla u\nabla\phi_jdx=a\int_{\Omega}u\phi_jdx+b\int_{\Omega}v\phi_jdx$$
	and
	$$\lambda_j\int_{\Omega}v\phi_jdx=\int_{\Omega}\nabla v\nabla\phi_jdx=c\int_{\Omega}u\phi_jdx+d\int_{\Omega}v\phi_jdx,$$
	in the matricial form,
	$$Az=\lambda_jz,\mbox{ where }z=\left(\begin{array}{c}
	\int_{\Omega}u\phi_jdx\\	\int_{\Omega}v\phi_jdx\end{array}\right).
	$$
	That is, $\lambda_j$ is eigenvalue of $A$ with eigenvector (nonzero) $z=\left(\begin{array}{c}
	\int_{\Omega}u\phi_jdx\\	\int_{\Omega}v\phi_jdx\end{array}\right)$.
	
	Therefore, $A$ has at most two eigenvalues of $(-\Delta,H_{0}^{1}(\Omega))$. If $A$ has two eigenvalues of $(-\Delta,H_{0}^{1}(\Omega))$, $\lambda_j$ and $\lambda_m$, we have
	$$\int_{\Omega}u\phi_kdx=0,\mbox{ for all }k\neq j,m.$$
	We conclude that there are, in the maximum, two eigenvalues, consider $\lambda_j$ and $\lambda_m$ such that $u=\alpha_1\phi_j+\beta_1\phi_m$ and $v=\alpha_2\phi_j+\beta_2\phi_m$, where
	$$\left(\begin{array}{c}
	\alpha_1\\	\alpha_2 \end{array}\right)\mbox{ is eigenvector of $A$ associated with the eigenvalue }\lambda_j$$
	and
	$$\left(\begin{array}{c}
	\beta_1\\	\beta_2 \end{array}\right)\mbox{is eigenvector of $A$ associated with the eigenvalue }\lambda_m.$$
	The itens $(ii)$ and $(iii)$, follow analyzing what has been done up.
	
	All details can be found in \cite{tese}, doctoral thesis that deals with this subject in detail.
\end{dem}	
	
Now, suppose that $(Q_1)$ admits a nonegative and nonzero solution, in the sense $U\geq0$ if $u\geq0$ and $v\geq0$. By demonstration of the above lemma, we have $u=\alpha_1\phi_j+\beta_1\phi_m$ and $v=\alpha_2\phi_j+\beta_2\phi_m$. We claim that, $j=1$ or $m=1$. Indeed, otherwise, as $\int_{\Omega}u\phi_1dx>0$ and from orthogonality of the eigenfunctions associated with distinct eigenvalues of $(-\Delta,H_{0}^{1}(\Omega))$, follow that
$$0<\int_{\Omega}u\phi_1dx=\alpha_1\int_{\Omega}\phi_j\phi_1dx+\beta_1\int_{\Omega}\phi_m\phi_1dx=0$$
which is a absurd.

Suppose that $j=1$. Thus, $u=\alpha_1\phi_1+\beta_1\phi_m$ and $v=\alpha_2\phi_1+\beta_2\phi_m$, so
$$0<\int_{\Omega}u\phi_1dx=\alpha_1\int_{\Omega}\phi_1^2dx$$
that is, $\alpha_1>0$. Analogously, $\alpha_2>0$. And, consequently, there is an eigenvector of $A$, associated with $\lambda_1$, having both positive coordinates.

From these comments, we have the proof of the Lemma \ref{EE}. Futhermore, recalling that $\frac{\partial\phi_1}{\partial\eta}<0$ on $\partial\Omega$, we have the proof of the Corollary \ref{CP}.

\end{document}